\documentclass[12pt]{article}
\usepackage{e-jc}
\dateline{Aug 16, 2019}{Jan 31, 2021}{Feb 12, 2021}

\usepackage{amsmath,amssymb,amsthm}

\usepackage{cleveref}
\crefname{equation}{}{}
\Crefname{equation}{Equation}{Equations}
\crefname{theorem}{Theorem}{Theorems}
\Crefname{theorem}{Theorem}{Theorems}
\crefname{lemma}{Lemma}{Lemmas}
\Crefname{lemma}{Lemma}{Lemmas}
\crefname{proposition}{Proposition}{Propositions}
\Crefname{proposition}{Proposition}{Propositions}
\crefname{corollary}{Corollary}{Corollaries}
\Crefname{corollary}{Corollary}{Corollaries}
\crefname{conjecture}{Conjecture}{Conjectures}
\Crefname{conjecture}{Conjecture}{Conjectures}
\crefname{section}{Section}{Sections}
\Crefname{section}{Section}{Sections}
\crefname{example}{Example}{Examples}
\Crefname{example}{Example}{Examples}
\crefname{problem}{Problem}{Problems}
\Crefname{problem}{Problem}{Problems}
\crefname{remark}{Remark}{Remarks}
\Crefname{remark}{Remark}{Remarks}
\crefname{figure}{Figure}{Figures}
\Crefname{figure}{Figure}{Figures}

\newcommand{\norm}[1]{\left\lVert#1\right\rVert}

\numberwithin{equation}{section}
\numberwithin{theorem}{section}

\Copyright{The author. Released under the CC BY-ND license (International 4.0).}

\title{Inequalities for doubly nonnegative functions}

\author{
Alexander Sidorenko\\
\small R\'{e}nyi Institute\\[-0.8ex]
\small Budapest, Hungary\\
\small\tt sidorenko.ny@gmail.com
}

\begin{document}

\maketitle

\begin{abstract}
Let $g$ be a bounded symmetric measurable nonnegative 
function on $[0,1]^2$, 
and $\norm{g} = \int_{[0,1]^2} g(x,y) dx dy$. 
For a graph $G$ with vertices $\{v_1,v_2,\ldots,v_n\}$ 
and edge set $E(G)$, we define
\[
  t(G,g) \; = \;
  \int_{[0,1]^n} \prod_{\{v_i,v_j\} \in E(G)} g(x_i,x_j)
    \: dx_1 dx_2 \cdots dx_n \; .
\]
We conjecture that $t(G,g) \geq \norm{g}^{|E(G)|}$ 
holds for any graph $G$ 
and any function $g$ with nonnegative spectrum.  
We prove this conjecture for various graphs $G$, 
including complete graphs, unicyclic and bicyclic graphs,
as well as graphs with $5$ vertices or less.
%
\end{abstract}

\section{Introduction}

Let $\mu$ be the Lebesgue measure on $[0,1]$.
Let $\mathcal{H}$ 
denote the space of bounded measurable real functions on $[0,1]^2$, 
and $\mathcal{G} \subset \mathcal{H}$ 
denote the subspace of symmetric functions. 
Let $\mathcal{H}_+$ and $\mathcal{G}_+$ 
denote the subsets of nonnegative functions in $\mathcal{H}$ 
and $\mathcal{G}$, respectively. 

Let $G$ be a simple graph with vertices $\{v_1,v_2,\ldots,v_n\}$ 
and edge set $E(G)$. 
We would like to know what conditions on $G$ and $g\in\mathcal{G}_+$ 
guarantee that 
\begin{equation}\label{eq:CI2}
  t(G,g) \overset{\underset{\mathrm{def}}{}}{=}
  \int_{[0,1]^n} \prod_{\{v_i,v_j\} \in E(G)} g(x_i,x_j) 
    \: d\mu^n
      \; \geq \;
  \left( \int_{[0,1]^2} g \: d\mu^2 \right)^{|E(G)|} .
\end{equation}
One approach is to ask what graphs $G$ satisfy \cref{eq:CI2} 
for every function $g\in\mathcal{G}_+$. 
It is easy to show that such graphs can not have odd cycles, 
so only graphs with chromatic number $2$ are suitable candidates. 
It led to 
\begin{conjecture}[\cite{Sidorenko:1992,Sidorenko:1993}]
\label{conj:SC}
Let $H$ be a bipartite graph with two vertex sets 
$V=\{v_1,v_2,\linebreak\ldots,v_n\}$, $W=\{w_1,w_2,\ldots,w_m\}$ 
and edge set $E(H) \subseteq V \times W$. 
Then for any function $h\in\mathcal{H}_+$
(not necessarily symmetric)
\begin{equation}\label{eq:SC}
  t(H,h) \overset{\underset{\mathrm{def}}{}}{=}
  \int_{[0,1]^{n+m}} \prod_{(v_i,w_j) \in E(H)} h(x_i,y_j) 
    \: d\mu^{n+m}
      \; \geq \;
  \left( \int_{[0,1]^2} h \: d\mu^2 \right)^{|E(H)|} .
\end{equation}
\end{conjecture}
We discuss \cref{conj:SC} in \cref{sec:SC}. 

\vspace{2mm}
For a (simple or bipartite) graph $G$, 
let $E(G)$ denote its edge set, and ${\rm e}(G)=|E(G)|$. 
For a simple graph $G$, 
let $V(G)$ denote its vertex set, and ${\rm v}(G)=|V(G)|$.

The 1-{\it subdivision} of a simple graph $G$ 
is a bipartite graph $H={\rm Sub}(G)$ 
with vertex sets $V(G)$ and $E(G)$, 
where $v \in V(G)$ and $e \in E(G)$ form an edge in $H$ 
if $v \in e$ in $G$.

We call a bipartite graph $H$ {\it symmetric} 
if it has an automorphism $\phi$ 
which switches its vertex-sets $V$ and $W$:\; 
$\phi(V)=W$, $\phi(W)=V$. 

We call a function $g\in\mathcal{G}_+$ {\it doubly nonnegative} 
if there is a function $h\in\mathcal{H}$ such that 
$g(x,y) = \int_{[0,1]} h(x,z)h(y,z) d\mu(z)$. 
Equivalently, a doubly nonnegative function 
is a nonnegative symmetric function with nonnegative spectrum. 
We call a function $g\in\mathcal{G}_+$ {\it completely positive} 
if there is a function $h\in\mathcal{H}_+$ such that 
$g(x,y) = \int_{[0,1]} h(x,z)h(y,z) d\mu(z)$. 
The terms ``doubly nonnegative'' and ``completely positive'' come 
from matrix theory; 
there exist functions which are doubly nonnegative 
but not completely positive 
(see \cref{sec:matrices}). 

In this article, we study two problems: 
(a) what functions $g\in\mathcal{G}_+$ satisfy 
$t(G,g) \geq \norm{g}^{{\rm e}(G)}$ 
for all simple graphs $G$ (we call such functions {\it nice}); and 
(b) what graphs $G$ satisfy the same inequality 
for any doubly nonnegative function $g$ 
(we call such graphs {\it good}).

If $G$ is good, then 
\cref{conj:SC} holds for $H={\rm Sub}(G)$. 
We show in \cref{sec:SC} that for a fixed $G$, 
inequality (\ref{eq:CI2}) 
holds for any completely positive function $g$ 
if and only if 
\cref{conj:SC} holds for $H={\rm Sub}(G)$. 
Thus, it is reasonable to expect that 
all completely positive functions are nice.

\begin{conjecture}\label{conj:CI}
All doubly nonnegative functions are nice. 
All simple graphs are good.
\end{conjecture}

Our \cref{th:permutation} demonstrates 
that there are nice functions which are not doubly nonnegative. 

If chromatic number $\chi(G)=2$, 
then goodness of $G$ should follow from \cref{conj:SC}.
In Sections \ref{sec:norming}-\ref{sec:small}, 
we give examples of good graphs $G$ 
with $\chi(G)\geq 3$. 
In particular, we prove that complete graphs, 
graphs whose complements consist of disjoint edges,
unicyclic and bicyclic graphs, 
generalized theta graphs,
and graphs with $\leq 5$ vertices 
are all good.  

In \cref{sec:extra-good}, 
we consider a strengthened variant of 
inequality~(\ref{eq:CI2}) which in many instances 
is easier to prove than the original one.
We say that a simple graph $G$ 
is {\it extra-good} if for any doubly nonnegative function $g$ 
and any bounded measurable nonnegative functions 
$f_1,\ldots,f_{{\rm v}(G)}$ on $[0,1]$,
\begin{equation}\label{eq:extra-good}
  \int_{[0,1]^{{\rm v}(G)}} \prod_{\{v_i,v_j\} \in E(G)} g(x_i,x_j) \: 
  \prod_{i=1}^{{\rm v}(G)} f_i(x_i)
    \: d\mu^{{\rm v}(G)}
      \geq 
  \left( \int_{[0,1]^2} f(x) g(x,y) f(y) 
    \: d\mu^2 \right)^{{\rm e}(G)} ,
\end{equation}
where 
$f(x) = \left(\prod_{i=1}^{{\rm v}(G)} 
  f_i(x)\right)^{1/(2{\rm e}(G))}$. 
Obviously, if $G$ is extra-good, then $G$ is good. 
For each graph that we proved to be good, 
we also were able to prove that it is extra-good.
It is possible that every good graph is extra-good.

A connection to Kohayakawa--Nagle--R\"{o}dl-Schacht conjecture 
is discussed in \cref{sec:KNRS}. 

In \cref{sec:multivar}, we discuss generalizations 
of inequality \cref{eq:CI2} 
for bounded measurable nonnegative symmetric functions 
of $r \geq 3$ variables.

\section{Doubly nonnegative and completely positive matrices}
\label{sec:matrices}

A {\it doubly nonnegative matrix} 
is a real positive semidefinite square matrix 
with nonnegative entries. 
A {\it completely positive matrix} 
is a doubly nonnegative matrix which can be factorized as $A=BB^{T}$ 
where $B$ is a nonnegative (not necessarily square) matrix. 
It is well known (see \cite{Berman:2003}) 
that for any $k\geq 5$ 
there exist doubly nonnegative $k \times k$ matrices 
which are not completely positive.

For a $k \times k$ matrix $A=[a_{ij}]$, 
we define a function $g_A$ on $[0,1]^2$ as 
$g(x,y)=a_{ij}$ for $(i-1)/k < x \leq i/k$, 
$(j-1)/k < y \leq j/k$, 
and $g(x,y)=0$ if $xy=0$. 
Obviously, 
$g_A$ is a doubly nonnegative 
(completely positive) function
if and only if 
$A$ is a doubly nonnegative 
(completely positive) matrix. 

If $A(H)$ is the adjacency matrix 
of a $k$-vertex simple graph $H$,
then the number of homomorphisms from $G$ to $H$ is equal to 
$t(G,g_{A(H)}) k^{{\rm v}(G)}$. 

Notice that if a nonzero $k \times k$ matrix $A$ 
has zero diagonal, 
then $g_A$ is not nice, 
since $t(G,g_A)=0$ for any graph $G$ with chromatic number $\chi(G)>k$. 
We are going to demonstrate that 
presence of a single positive diagonal entry 
can be sufficient to make $g_A$ nice.

\begin{theorem}\label{th:permutation}
Let $P$ be a symmetric permutation matrix of order $k$ 
with $a \geq 1$ diagonal entries equal to $1$, 
and $b \geq 1$ pairs of off-diagonal entries equal to $1$ $(a+2b=k)$. 
Then $g_P$, while not being positive semidefinite, is a nice function. 
\end{theorem}

\begin{proof}[\bf{Proof}]
$P$ has eigenvalues $1$ with multiplicity $a+b$, 
and $-1$ with multiplicity $b\geq 1$. 
Therefore, $P$ is not positive semidefinite.
If graph $G$ has connected components $G_1,G_2,\ldots,G_m$ 
then $t(G,g) = \prod_{i=1}^m t(G_i,g)$. 
Hence, to prove \cref{eq:CI2} for $g=g_P$ 
it is sufficient to consider connected graphs $G$. 
If $G$ is a tree, 
then validity of \cref{eq:CI2} follows from \cref{eq:SC} 
(\cref{conj:SC} has been proved for trees by various authors; 
a short proof can be found in \cite{Jagger:1996}). 
Hence, we may assume that $G$ is not a tree. 
If $n={\rm v}(G)$, then 
${\rm e}(G) \geq n$. 
As $P$ has $a \geq 1$ diagonal entries equal to $1$, we get 
$t(G,g_P) \geq a \: (1/k)^n \geq (1/k)^{{\rm e}(G)}$ 
and $\int_{[0,1]^2} g_P \: d\mu^2 = 1/k$.
\end{proof}

\section{More on \cref{conj:SC}}\label{sec:SC}

The earliest known works where inequalities of type 
\cref{eq:CI2,eq:SC} appear 
are \cite{Mulholland:1959} and \cite{Atkinson:1960}. 
In 1959, Mulholland and Smith \cite{Mulholland:1959} 
proved that for any symmetric nonnegative matrix $A$ 
and any nonnegative vector ${\bf z}$,
\begin{equation}\label{eq:MS}
  ({\bf z}^{\mathsf{T}}\! A^k {\bf z}) \cdot 
  ({\bf z}^{\mathsf{T}} {\bf z})^{k-1}
    \; \geq \; 
  ({\bf z}^{\mathsf{T}}\! A {\bf z})^k \; ,
\end{equation}
where equality takes place if and only if 
${\bf z}$ is an eigenvector of $A$ or a zero vector.
Note that (\ref{eq:MS}) 
is a particular case of (\ref{eq:CI2}) 
where $H$ is the $k$-edge path $P_k$. 

Almost at the same time, 
Atkinson, Watterson and Moran \cite{Atkinson:1960} 
proved that
$
  nm \cdot s(A A^{\mathsf{T}}\! A) \; \geq \; s(A)^3 \: ,
$
where $A$ is an (asymmetric) nonnegative $(n \times m)$-matrix, 
and $s(A)$ is the sum of entries of $A$. 
They presented their inequality in both matrix and integral form, 
and conjectured validity of \cref{eq:SC} 
for $H=P_k$ with $k \geq 3$. 

In 1965, Blakley and Roy \cite{Blakley:1965},
being unaware of the article \cite{Mulholland:1959}, 
rediscovered (\ref{eq:MS}). 

Lately, \cref{conj:SC} has been proved for various bipartite graphs 
(see \cite{Conlon:2010,Conlon:2018,Conlon:2017,Conlon:2019,Kim:2016,
Li:2011,Lovasz:2011,Parczyk:2014,Sidorenko:1992,Sidorenko:1993,
Sidorenko:1991,Szegedy:2014,Szegedy:2015}), 
among them: trees, complete bipartite graphs, 
and graphs with $9$ vertices or less. 
Some of the authors restricted \cref{eq:SC} to symmetric functions $h$. 
Nevertheless, the proofs of their results can be extended 
to asymmetric $h$ as well. 
Let $\mathfrak{S}$ be the class of bipartite graphs 
that satisfy \cref{conj:SC}, 
and $\mathfrak{S}_*$ be the class of bipartite graphs $H$ 
that satisfy \cref{eq:SC} for all $h\in\mathcal{G}_+$. 
Obviously, $\mathfrak{S}\subseteq\mathfrak{S}_*$. 
It would be nice to prove 
$\mathfrak{S} = \mathfrak{S}_*$. 

\begin{theorem}\label{th:symm}
If $H\in\mathfrak{S}_*$ is symmetric, then $H\in\mathfrak{S}$.
\end{theorem}

In the proof of \cref{th:symm}, we will use 
the so called ``tensor-trick'' lemma.

\begin{lemma}[\cite{Sidorenko:1992,Sidorenko:1993}]
\label{th:tensor_trick}
If there exists a constant $c = c_H > 0$ such that 
for any $h\in\mathcal{H}_+\,$, 
$t(H,h) \geq c \cdot \left(\int_{[0,1]^2} h d\mu^2\right)^{{\rm e}(H)}$, 
then $H\in\mathfrak{S}$.
\end{lemma}

\begin{proof}[\bf{Proof of \cref{th:symm}}]
It is sufficient to consider the case when $H$ is connected. 
Denote by $n$ the size of each vertex set of $H$, 
so the total number of vertices is $2n$. 
Let $h\in\mathcal{H}_+$. 
Define its ``transpose'' $h^{\mathsf{T}}$ as 
$h^{\mathsf{T}}(x,y) = h(y,x)$. 
As $H$ is symmetric, $t(H,h)=t(H,h^{\mathsf{T}})$. 
Define symmetric function $\tilde{h}\in\mathcal{G}_+$ as follows:
\[
  \tilde{h}(x,y) = 
  \begin{cases} 
             0 & \mbox{if } 0 \leq x,y < 1/2; \\
    h(2x,2y-1) & \mbox{if } 0 \leq x < 1/2 \leq y \leq 1; \\
    h(2y,2x-1) & \mbox{if } 0 \leq y < 1/2 \leq x \leq 1; \\
             0 & \mbox{if } 1/2 \leq x,y \leq 1.
  \end{cases}
\]
Notice that 
$\int_{[0,1]^2} \tilde{h} d\mu^2 = (1/2) \int h_{[0,1]^2} \: d\mu^2$. 
As $H$ is connected, 
\[
  t(H,\tilde{h})
  = 2^{-2n} (t(H,h) + t(H,h^{\mathsf{T}})) 
  = 2^{1-2n} \: t(H,h) .
\]
Since $H \in \mathfrak{S}_{*}$, we get
$t(H,\tilde{h}) \geq \left(\int_{[0,1]^2} \tilde{h} d\mu^2\right)^{{\rm e}(H)}$.
Hence, 
\[
  t(H,h) \geq 
  2^{2n-1-{\rm e}(H)} 
    \left(\int_{[0,1]^2} h d\mu^2\right)^{{\rm e}(H)} ,
\]
and by \cref{th:tensor_trick},  $H\in\mathfrak{S}$.
\end{proof}

\begin{remark}\label{th:measure}
It is a classical fact that there exists a measure preserving bijection 
between any two atomless measure spaces with total measure $1$. 
In particular, if $\mu_1$ and $\mu_2$ are atomless measures on $[0,1]$, 
and a bipartite graph $H\in\mathfrak{S}$ 
has vertex sets of sizes $n$ and $m$, 
then for any bounded non-negative function $h$ on $[0,1]^2$, 
measurable with respect to $\mu_1 \otimes \mu_2$,
\[
  \int_{[0,1]^{n+m}} \prod_{(v_i,w_j) \in E(H)} h(x_i,y_j) 
    \: d\mu_1^n d\mu_2^m
      \; \geq \;
  \left( \int_{[0,1]^2} h \: d\mu_1 d\mu_2 \right)^{{\rm e}(H)} .
\]
\end{remark}

\begin{theorem}\label{th:same_degree}
If \cref{conj:SC} holds for a bipartite graph $H$, 
and all vertices from the first vertex set of $H$ 
have the same degree $a$, 
then $t(H,h) \geq t(K_{1,a},h)^{{\rm e}(H)/a}$ for any $h\in\mathcal{H}_+$. 
If all vertices from the second vertex set of $H$ 
have the same degree $b$, 
then $t(H,h) \geq t(K_{b,1},h)^{{\rm e}(H)/b}$ for any $h\in\mathcal{H}_+$. 
\end{theorem}

\begin{proof}[\bf{Proof}]
We will prove the first part of the statement 
(the proof of the second part is similar). 
Notice that ${\rm e}(H)=na$. 
It is sufficient to consider functions 
$h\in\mathcal{H}_+$ that are separated from zero: 
$\inf_{[0,1]^2} h > 0$. 
Denote $\varphi(x)=\int_{[0,1]} h(x,y)d\mu(y)$. 
Then $c = \int_{[0,1]} \varphi(x)^{a} d\mu > 0$,
and $f(x) = \varphi(x)^{a} / c$ 
is positive and bounded on $[0,1]$. 
Consider a measure $\mu_*$ on $[0,1]$ 
defined by $d\mu_* = fd\mu$, 
so $\mu_*([0,1])=1$. 
Denote $\widehat{h}(x,y) = h(x,y) f(x)^{-1/a}$. 
Clearly, $\widehat{h}$ is bounded and measurable 
with respect to $\mu_* \otimes \mu$. 
By~\cref{th:measure},
\begin{eqnarray*}
  t(H,h)^{1/n} 
  & = &
  \left( \int_{[0,1]^{n+m}}
    \prod_{(v_i,w_j) \in E(H)} h(x_i,y_j) 
    d\mu^n d\mu^m
  \right)^{1/n}
  \\
  & = &
  \left( \int_{[0,1]^{n+m}}
    \prod_{(v_i,w_j) \in E(H)} \widehat{h}(x_i,y_j) 
    d\mu_*^n d\mu^m
  \right)^{1/n}
  \\
  & \geq &
  \left( \int_{[0,1]^2} \widehat{h}(x,y) d\mu_*(x) d\mu(y) \right)^a
  \\
  & = &
  \left( \int_{[0,1]^2} h(x,y) f(x)^{-1/a} d\mu_*(x) d\mu(y) \right)^a
  \\
  & = &
  \left( \int_{[0,1]^2} h(x,y) f(x)^{1-(1/a)} d\mu^2 \right)^a
  \\
  & = &
  \left( \int_{[0,1]} \varphi(x) f(x)^{1-(1/a)} d\mu \right)^a
  \\
  & = &
  c^{1-a} \left( \int_{[0,1]} \varphi(x)^a d\mu \right)^a
  \; = \;
  \int_{[0,1]} \varphi(x)^a d\mu 
  \; = \;
  t\left(K_{1,a},h\right)
  \; .
\end{eqnarray*}
\end{proof}

\begin{theorem}\label{th:sub}
For a fixed graph $G$, 
inequality \cref{eq:CI2} holds for any completely positive function $g$ 
if and only if \cref{conj:SC} holds for $H={\rm Sub}(G)$. 
\end{theorem}

\begin{proof}[\bf{Proof}]
Suppose that \cref{conj:SC} holds for $H={\rm Sub}(G)$, 
and a function $g$ is completely positive. 
There exists $h\in\mathcal{H}_+$ such that 
$g(x,y) = \int_{[0,1]} h(x,z)h(y,z) d\mu(z)$. 
Then $t(G,g)=t(H,h)$. 
Every vertex in the second vertex set of $H$ has degree $2$. 
By \cref{th:same_degree}, we have $t(H,h) \geq t(K_{2,1},h)^{{\rm e}(H)/2}$. 
As ${\rm e}(G) = {\rm e}(H)/2$ and $t(K_{2,1},h) = \int_{[0,1]^2} g d\mu^2$, 
we get \cref{eq:CI2}.

Now suppose that \cref{eq:CI2} holds 
for any completely positive function $g$. 
Let $h\in\mathcal{H}_+$. 
Set $g(x,y) = \int_{[0,1]} h(x,z)h(y,z) d\mu(z)$. 
Then $t(H,h) = t(G,g) \geq (\int_{[0,1]^2} g d\mu^2)^{{\rm e}(G)}
 \geq (\int_{[0,1]^2} h d\mu^2)^{2{\rm e}(G)}
 = (\int_{[0,1]^2} h d\mu^2)^{{\rm e}(H)}$.
\end{proof}

\section{Subdivisions that are norming}\label{sec:norming}

We say that a bipartite graph $H$ with vertex sets 
$V=\{v_1,v_2,\ldots,v_n\}$ and $W=\{w_1,w_2,\linebreak\ldots,w_m\}$ 
has the {\it H\"{o}lder property} if for any assignment 
$f: E(H) \to \mathcal{H}$,
\begin{equation}\label{eq:norming}
  \left(\int_{[0,1]^{n+m}} \prod_{e = (v_i,w_j) \in E(H)}
    f_e(x_i,y_j) 
    \: d\mu^{n+m} \right)^{{\rm e}(H)}
      \leq 
  \prod_{e \in E(H)} t(H, f_e) .
\end{equation}
It is known (see \cite{Hatami:2010,Lovasz:2010}) 
that every graph $H$ with the H\"{o}lder property 
(except a star with even number of edges)
is a {\it norming graph}: 
$t(H,h)^{1/{\rm e}(H)}$ is a norm on $\mathcal{H}$. 
Conversely, every norming graph has the H\"{o}lder property.

\begin{theorem}\label{th:norming}
If ${\rm Sub}(G)$ has the H\"{o}lder property then $G$ is good.
\end{theorem}

\begin{proof}[\bf{Proof}]
Let $H={\rm Sub}(G)$. 
If $g(x,y) = \int_{[0,1]} h(x,z)h(y,z) d\mu(z)$ 
then $t(G,g) = t(H,h)$. 
Select a pair of edges $e',e''$ in $H$ 
which subdivide the same edge of $G$. 
Assign $f_{e'} = f_{e''} = h$ and $f_e = 1$ for $e \neq e',e''$. 
Then the left hand side of \cref{eq:norming} is 
$\left(\int_{[0,1]^2} g d\mu^2\right)^{2{\rm e}(G)}$, 
and the right hand side is $t(G,g)^2$.
\end{proof}

The 1-subdivision of cycle $C_n$ is an even cycle $C_{2n}$ 
which is a norming graph. 
The 1-subdivision of the octahedron $K_{2,2,2}$ is norming 
(see \cite[Example 4.15]{Conlon:2017}). 
Hence, $C_n$ and $K_{2,2,2}$ are good graphs.

In a norming graph, the degrees of vertices are even 
(see Observation~2.5 in \cite{Hatami:2010}).
Hence, ${\rm Sub}(K_{2r})$ is not norming. 
While ${\rm Sub}(K_3)=C_6$ is norming, 
it is not known whether ${\rm Sub}(K_{2r+1})$ with 
$r \geq 2$ is norming. 
We will prove in the next section
that all complete graphs are good.

\section{Complete graphs are good}\label{sec:clique}

\begin{theorem}\label{th:G+v}
If graph $G$ is good, 
and graph $G_1$ is obtained from $G$ by adding a new vertex 
adjacent to all vertices of $G$, then $G_1$ is good.
\end{theorem}

\begin{corollary}\label{th:clique}
The complete graphs are good.
\end{corollary}

\begin{theorem}\label{th:G+vw}
If graph $G$ is good, 
and graph $G_2$ is obtained from $G$ by adding two vertices 
adjacent to all vertices of $G$ but not to each other, 
then $G_2$ is good.
\end{theorem}

\begin{corollary}
Any graph whose complement is a set of independent edges 
is good.
\end{corollary}

To prove \cref{th:G+v,th:G+vw} 
we need a couple of auxiliary results. 
For $g\in\mathcal{G}_+$, 
we denote $\norm{g} = \int_{[0,1]^2} g d\mu^2$. 

\begin{lemma}\label{th:m_1_2}
If function $g$ is doubly nonnegative, then 
$t(K_3,g) \! \geq t(K_{1,2},g)^{3/2}$
and \linebreak
$t(K_4-e,g) \geq t(K_{2,2},g)^{5/4}$.
\end{lemma}

\begin{proof}[\bf{Proof}]
As $g$ is doubly nonnegative, 
$g(x,y) = \int_{[0,1]} h(x,z)h(y,z) d\mu(z)$. 
Then $t(K_3,g)=t(C_6,h)$, $t(K_{1,2},g)=t(P_4,h)$, 
$t(K_{2,2},g)=t(C_8,h)$, 
where $P_4$ denotes the $4$-edge path.

As $C_6$ is norming, and $P_4$ is a subgraph of $C_6$, 
we have $t(C_6,h)^{1/6} \geq t(P_4,h)^{1/4}$.

By the Cauchy--Schwarz inequality, 
$t(K_4-e,g) \norm{g} \geq t(K_3,g)^2 \!= t(C_6,h)^2$\!\!. 
As $t(C_{2k},h)^{1/2k}$ is the $(2k)$-th Schatten norm of $h$, 
we get $t(C_6,h)^{1/6} \geq t(C_8,h)^{1/8}$. 
Hence, 
$t(K_4-e,g) \geq t(K_{2,2},g)^{3/2} \norm{g}^{-1}
 \geq t(K_{2,2},g)^{5/4}$.
\end{proof}

For $n > m \geq 1$, let $K_n-K_m$ 
denote the complement of $K_m$ in $K_n$. 
It follows from \cref{th:m_1_2} that 
for $m=1,2$ and any doubly nonnegative function $g$,
\begin{equation}\label{eq:m_1_2}
  t(K_{m+2}-K_m,g) \; \geq \; t(K_{2,m},g) \cdot \norm{g} \; .
\end{equation}
Thus, \cref{th:G+v,th:G+vw} follow from the next proposition.

\begin{proposition}\label{th:G_m}
Let integer $m$ be such that \cref{eq:m_1_2} 
holds for all doubly nonnegative functions. 
If graph $G$ is good, 
and graph $G_m$ is obtained from $G$ by adding 
a group of $m$ independent vertices that are 
adjacent to all vertices of $G$, 
then $G_m$ is good.
\end{proposition}

\begin{proof}[\bf{Proof}]
We may assume ${\rm v}(G) \geq 2$. 
It is sufficient to consider functions 
$g$ that are separated from zero: 
$\inf_{[0,1]^2} g > 0$. 
Then function 
$\varphi(x_1,\ldots,x_m) = 
 \int_{[0,1]} \prod_{i=1}^m g(x_i,y) d\mu(y)$ 
is positive and bounded on $[0,1]^m$. 
For each ${\bf x}=(x_1,\ldots,x_m) \in [0,1]^m$, 
consider measure $\mu_{\bf x}$ on $[0,1]$, 
defined by $d\mu_{\bf x} = f_{\bf x} d\mu$, 
where 
$f_{\bf x}(y) =
 \prod_{i=1}^m g(x_i,y) \varphi(x_1,\ldots,x_m)^{-1}$. 
It is easy to see that $\mu_{\bf x}([0,1])=1$, 
and $g$ is bounded and measurable 
with respect to $\mu_{\bf x} \otimes \mu_{\bf x}$. 
By~\cref{th:measure}, 
\[
  t(G,g,\mu_{\bf x}) \overset{\underset{\mathrm{def}}{}}{=} 
  \int_{[0,1]^{{\rm v}(G)}}
    \prod_{\{v_i,v_j\} \in E(G)} g(y_i,y_j) 
    \: d\mu_{\bf x}^{{\rm v}(G)}
      \; \geq \;
  \left( \int_{[0,1]^2} g \: d\mu_{\bf x}^2 \right)^{{\rm e}(G)}
  .
\]
Hence,
\begin{multline*}
  t(G_m,g) =  
  \int_{[0,1]^m} t(G,g,\mu_{x_1,\ldots,x_m}) 
    \varphi(x_1,\ldots,x_m)^{{\rm v}(G)} 
      \: d\mu(x_1) \cdots d\mu(x_m)
  \\  \geq  
  \int_{[0,1]^m} 
    \left( \int_{[0,1]^2} g \: d\mu_{x_1,\ldots,x_m}^2
    \right)^{{\rm e}(G)} 
    \varphi(x_1,\ldots,x_m)^{{\rm v}(G)} 
      \: d\mu(x_1) \cdots d\mu(x_m)
  \\  = 
  \int_{[0,1]^m} 
    \left( \int_{[0,1]^2} 
      g(y,z) \prod_{i=1}^m \left( g(x_i,y) g(x_i,z) \right)
      \: d\mu(y) d\mu(z) \right)^{{\rm e}(G)} 
  \\ \times
    \varphi(x_1,\ldots,x_m)^{{\rm v}(G) - 2{\rm e}(G)} 
      \: d\mu(x_1) \cdots d\mu(x_m)
    .
\end{multline*}
As $1 + {\rm e}(G) - {\rm v}(G)/2 \leq {\rm e}(G)$, 
by using the H\"{o}lder inequality, we get
\begin{align*}
  t(G_m,g) \cdot & t(K_{2,m},g)^{{\rm e}(G) - {\rm v}(G)/2} 
  \\ & = \;
  t(G_m,g) \cdot 
  \left(\int_{[0,1]^m} 
    \varphi^2(x_1,\ldots,x_m) 
      \: d\mu(x_1) \cdots d\mu(x_m)
  \right)^{{\rm e}(G) - {\rm v}(G)/2}
  \\ & \geq \; 
  \left(\int_{[0,1]^{m+2}} \!
      g(y,z) \prod_{i=1}^m \left( g(x_i,y) g(x_i,z) \right)
      \: d\mu(y) d\mu(z) 
      d\mu(x_1) \cdots d\mu(x_m)
  \right)^{{\rm e}(G)} 
  \\ & = \;
  t(K_{m+2}-K_m,g)^{{\rm e}(G)} \: .
\end{align*}
As \cref{conj:SC} holds for complete bipartite graphs
(see \cite{Sidorenko:1992}), $K_{2,m}$ is good. 
By using \cref{eq:m_1_2}, we get
\begin{align*}
  t(G_m,g) \: & \geq \: 
  t(K_{2,m},g)^{-{\rm e}(G) + {\rm v}(G)/2} \:
  t(K_{2,m},g)^{{\rm e}(G)} \norm{g}^{{\rm e}(G)}
    \\ & = \;
  t(K_{2,m},g)^{{\rm v}(G)/2} \norm{g}^{{\rm e}(G)} \; \geq \;
  \norm{g}^{m{\rm v}(G)} \norm{g}^{{\rm e}(G)} \; = \;
  \norm{g}^{{\rm e}(G_m)} .
\end{align*}
\end{proof}

\section{Extra-good graphs}\label{sec:extra-good}

\begin{theorem}
If $G$ is vertex-transitive and good, then $G$ is extra-good.
\end{theorem}

\begin{proof}[\bf{Proof}]
We denote $n={\rm v}(G)$, 
$f(x) = \left(\prod_{i=1}^n 
  f_i(x)\right)^{1/(2{\rm e}(G))}$ 
and 
$\tilde{g}(x,y) = 
f(x)g(x,y)f(y)$. 
Notice that $\tilde{g}$ is doubly nonnegative. 
Since $G$ is vertex-transitive, 
any permutation of the functions 
$f_1,f_2,\ldots,f_n$ 
does not change the value of the integral
on the left hand side of \cref{eq:extra-good}. 
By applying the H\"{o}lder inequality 
to the geometric mean of all $n!$ possible integrals, 
we get
\begin{align*}
  \int_{[0,1]^n} \prod_{\{v_i,v_j\} \in E(G)} g(x_i,x_j) \: 
  \prod_{i=1}^n f_i(x_i)
    \: d\mu^n
  \: & \geq 
  \int_{[0,1]^n} \prod_{\{v_i,v_j\} \in E(G)} g(x_i,x_j) \: 
  \prod_{i=1}^n f(x_i)^{2{\rm e}(G)/n}
    \: d\mu^n
  \\ & = 
  \int_{[0,1]^n} \prod_{\{v_i,v_j\} \in E(G)} \tilde{g}(x_i,x_j)
    \: d\mu^n
    \; .
\end{align*}
Since $G$ is good, 
\[
  \int_{[0,1]^n} \prod_{\{v_i,v_j\} \in E(G)} \tilde{g}(x_i,x_j)
    \: d\mu^n
    \; \geq \;
  \left( \int_{[0,1]^2} \tilde{g}(x,y) 
    \: d\mu^2 \right)^{{\rm e}(G)} .
\]
\end{proof}

\begin{corollary}
Cycles and complete graphs are extra-good.
\end{corollary}

\begin{theorem}\label{th:leaf}
If vertex $v$ is a leaf in graph $G$, 
and $G-v$ is extra-good, then $G$ is extra-good too.
\end{theorem}

\begin{proof}[\bf{Proof}]
Let ${\rm v}(G)=n+1$, $V(G)=\{v_1,\ldots,v_n,v_{n+1}\}$, 
and $v=v_{n+1}$ is a leaf. 
For a set of functions $f_1,\ldots,f_n,f_{n+1}$, 
define $\tilde{f_i} = f_i$ for $i=1,\ldots,n-1$, 
and 
$\tilde{f_n}(x) = 
f_n(x) \int_{[0,1]} g(x,y) \, f_{n+1}(y) \, d\mu(y)$. 
Then
\begin{align*}
  \int_{[0,1]^{n+1}} 
    \prod_{\{v_i,v_j\} \in E(G)} g(x_i,x_j) \:
    \prod_{i=1}^{n+1} f_i(x_i) d\mu^{n+1}
  \: & =
  \int_{[0,1]^n} 
    \prod_{\{v_i,v_j\} \in E(G-v)} g(x_i,x_j) \:
    \prod_{i=1}^n \tilde{f}_i(x_i) d\mu^n
  \\ & \geq   
  \left( \int_{[0,1]^2} \tilde{f}(x) g(x,y) \tilde{f}(y)
    \: d\mu^2 \right)^{{\rm e}(G)-1} 
  \\ & = \; I_0^{{\rm e}(G)-1},
\end{align*}
where
$\tilde{f}(x)=\left(\prod_{i=1}^n 
  \tilde{f}_i(x)\right)^{1/(2{\rm e}(G)-2)}$. 
Set $\varphi(x) = \int_{[0,1]} g(x,y) f_{n+1}(y) d\mu(y)$, 
\[
  I_1 = \int_{[0,1]^2} \varphi(x)^{-1} g(x,y) f_{n+1}(y) d\mu^2,
  \;\;\;\;\;\;
  I_2 = \int_{[0,1]^2} f_{n+1}(x) g(x,y) \varphi(y)^{-1} d\mu^2.
\]
It is easy to see that
$I_1 = \int_{[0,1]} \varphi(x)^{-1} \varphi(x) d\mu = 1$, 
and similarly, $I_2 = 1$. 
By the H\"{o}lder inequality, 
\[
  I_0^{{\rm e}(G)-1} \cdot 1 \cdot 1
    \; = \;
  I_0^{{\rm e}(G)-1} \cdot I_1^{1/2} \cdot I_2^{1/2}
    \; \geq \;
  \left( \int_{[0,1]^2} f(x) \: g(x,y) \: f(x) \: d\mu^2
  \right)^{{\rm e}(G)} ,
\]
where  
\begin{equation*}
  f^{2{\rm e}(G)}
    = 
  (\tilde{f})^{2{\rm e}(G)-2} \: \varphi^{-1} f_{n+1}
    =
  \prod_{i=1}^n \tilde{f}_i \cdot \varphi^{-1} f_{n+1}
    =
  \prod_{i=1}^{n-1} f_i \cdot (f_n \varphi) 
    \cdot \varphi^{-1} f_{n+1}
    = 
  \prod_{i=1}^{n+1} f_i \; .
\end{equation*}
\end{proof}

A connected graph $G$ is called {\it unicyclic} 
if ${\rm e}(G) = {\rm v}(G)$, 
and {\it bicyclic} if ${\rm e}(G) = {\rm v}(G)+1$. 

\begin{corollary}\label{th:unicyclic}
Trees and unicyclic graphs are extra-good.
\end{corollary}

Let $P(k_1,k_2,\ldots,k_r)$ denote a graph 
consisting of two vertices 
joined by $r$ internally disjoint paths 
of length $k_1,k_2,\ldots,k_r$. 
Such a graph is called {\it theta graph} 
(see \cite{Blinco:2004}). 
If $r=2$, then $P(k_1,k_2)$ is simply a cycle of length $k_1+k_2$.

\begin{theorem}\label{th:multipath}
$P(k_1,k_2,\ldots,k_r)$ is extra-good.
\end{theorem}

Let $C(k_1,k_2,m)$ denote  a graph which consists of 
two cycles of length $k_1$ and $k_2$ 
connected by a path of length $m \geq 0$ 
(when $m=0$, the cycles share a vertex). 
Notice that $P(k_1,k_2,k_3)$ and $C(k_1,k_2,m)$ 
are the only bicyclic graphs without leaves.
In view of \cref{th:leaf,th:multipath}, 
in order to prove that all bicyclic graphs are extra-good, 
it is sufficient to show that $C(k_1,k_2,m)$ is extra-good. 
This will follow from the next result.

\begin{theorem}\label{th:cycle}
Let graph $G_0$ be formed 
by attaching a cycle to one of the vertices of graph $G$. 
If $G$ is extra-good, then $G_0$ is extra-good.
\end{theorem}

\begin{corollary}\label{th:bicyclic}
Bicyclic graphs are extra-good.
\end{corollary}

The proofs of \cref{th:multipath,th:cycle} 
are given in the appendix.

Consider a tree $T$ whose vertices are arbitrarily colored 
in black and white so that at least one vertex is black. 
Take $r$ disjoint copies of $T$ and glue together 
``sister'' black vertices from different copies. 
We call the resulting graph a {\it multitree}. 
For example $P(k,k,\ldots,k)$ is a multitree. 
The case of even cycle in \cref{th:cycle} 
is a particular case of the following statement.

\begin{theorem}\label{th:multitree}
Let graph $G_0$ be formed 
by gluing a black vertex of a multitree 
to one of the vertices of graph $G$. 
If $G$ is extra-good, then $G_0$ is extra-good.
\end{theorem}

\begin{theorem}\label{th:G+v_extra}
If graph $G$ is extra-good, 
and graph $G_1$ is obtained from $G$ by adding a new vertex 
adjacent to all vertices of $G$, then $G_1$ is extra-good.
\end{theorem}

\begin{theorem}\label{th:G+vw_extra}
If graph $G$ is extra-good, 
and graph $G_2$ is obtained from $G$ by adding two vertices 
adjacent to all vertices of $G$ but not to each other, 
then $G_2$ is extra-good.
\end{theorem}

Theorems \ref{th:multitree}, 
\ref{th:G+v_extra}, and \ref{th:G+vw_extra} 
are not used in the rest of the article. 
We omit their proofs 
as they are very similar to the proofs of 
Theorems \ref{th:cycle}, \ref{th:G+v}, and \ref{th:G+vw}.

\section{Graphs with small number of vertices}\label{sec:small}

\begin{theorem}\label{th:small}
All graphs with $5$ vertices or less are good. 
\end{theorem}

\begin{proof}[\bf{Proof}]
If graph $G$ has connected components $G_1,G_2,\ldots,G_m$ 
then $t(G,g) \! = \! \prod_{i=1}^m t(G_i,g)$. 
Hence, it is sufficient to consider connected graphs $G$ only. 
As \cref{conj:SC} has been proved for bipartite graphs 
with $9$ vertices or less, it is sufficient to consider 
connected simple graphs with at least one odd cycle. 
The results of \cref{sec:clique,sec:extra-good} cover 
all such graphs with $\leq 5$ vertices.
A table of all $5$-vertex graphs 
that do not have isolated vertices 
can be found in \cite[Figure~6]{Adams:2008}.
\end{proof}

Some $6$-vertex graphs are not covered by 
the results of \cref{sec:clique,sec:extra-good}, 
but for almost all of them we were able to prove 
that they are good by using the Cauchy--Schwarz 
and H\"{o}lder inequalities. 
The only graph 
with $6$ vertices 
that we were unable to prove to be good 
is the complement of the $5$-edge path.

\section{Locally dense graphs}\label{sec:KNRS}

A simple graph $H$ is called $(\varepsilon,d)$-{\it dense} 
if every subset $X \subseteq V(H)$ of size 
$|X| \geq \varepsilon |V(H)|$ spans at least 
$\frac{d}{2} |X|^2$ edges. 

\begin{conjecture}[\cite{Kohayakawa:2010}]\label{conj:KNRS}
For any graph $G$ and $\delta,d \in (0,1)$, 
there exists $\varepsilon = \varepsilon(\delta,d,G)$ 
such that there are at least 
$(d^{{\rm e}(G)} - \delta) {\rm v}(H)^{{\rm v}(G)}$ 
homomorphisms of $G$ into any sufficiently large 
$(\varepsilon,d)$-dense graph $H$.
\end{conjecture}

When $\chi(G)=2$, 
\cref{conj:KNRS} follows from \cref{conj:SC}. 
It is known that \cref{conj:KNRS} holds for 
complete graphs and multipartite complete graphs 
(see\cite{Kohayakawa:2010,Reiher:2014}). 
Christian Reiher \cite{Reiher:2014} proved \cref{conj:KNRS} 
for odd cycles. 
Joonkyung Lee \cite{Lee:2019b} proved that
adding an edge to a cycle or a tree 
produces graphs that satisfy the conjecture. 
He also proved that \cref{conj:KNRS} holds 
for a class of graphs obtained by gluing 
complete multipartite graphs (or odd cycles) 
in a tree-like way. 
All graphs with $5$ vertices or less satisfy the conjecture.

For a nonnegative symmetric $k \times k$ matrix $A$, 
we define its {\it density} $d(A)$ as the minimum of 
${\bf x}^{\mathsf{T}} A {\bf x}$ over all nonegative $k$-dimensional 
vectors ${\bf x}$ with the sum of entries equal to $1$. 
Clearly, $d(A) \leq \norm{g_A}$. 
We call a graph $G$ {\it density-friendly} 
if $t(G,g_A) \geq d(A)^{{\rm e}(G)}$ for any $A$. 
It is easy to see that any density-friendly graph 
satisfies \cref{conj:KNRS}, but the converse is not obvious. 

While \cref{conj:CI,conj:KNRS} looks very different, 
the sets of graphs which are known to satisfy them 
are surprisingly similar.

One can try to build a bridge between these two topics 
by defining 
\[
  c(A) \: = \: \inf_G t(G,g_A)^{1/{\rm e}(G)} \: .
\]
Then \cref{conj:CI} claims that $c(A)=\norm{g_A}$ 
for any doubly nonnegative matrix $A$, 
and \cref{conj:KNRS} claims that $c(A) \geq d(A)$ 
for any nonnegative symmetric matrix $A$.

\section{Functions of $r$ variables}\label{sec:multivar}

Let $G$ be an $r$-uniform hypergraph with vertex set $V(G)=\{v_1,v_2,\ldots,v_n\}$ 
and edge set $E(G)$ 
(edges are $r$-element subsets of the vertex set). 
Let $g$ be a bounded symmetric measurable nonnegative 
function defined on $[0,1]^r$. 
Denote
$\norm{g} = \int_{[0,1]^r} g d\mu^r$ and 
\[
  t(G,g) \; = \;
  \int_{[0,1]^n} \prod_{\{v_{i_1},v_{i_2},\ldots,v_{i_r}\}\in E(G)}
    g(x_{i_1},x_{i_2},\ldots,x_{i_r})
    \: d \mu^n \; .
\]

\begin{problem}
Characterize functions $g$ such that 
\begin{equation}\label{eq:CIr}
  t(G,g) \; \geq \; \norm{g}^{|E(G)|} 
\end{equation}
holds for every $r$-uniform hypergraph $G$.
\end{problem}

When $r=1$, it is obvious that \cref{eq:CIr} 
holds for any nonnegative function $h$ on $[0,1]$. 

The {\it incidence graph} of an $r$-uniform hypergraph $G$ 
is a bipartite graph ${\rm Inc}(G)$ 
with vertex sets $V(G)$ and $E(G)$, 
where $v \in V(G)$ and $e \in E(G)$ 
form an edge $\{v,e\}$ in ${\rm Inc}(G)$ 
if and only if $v \in e$ in $G$. 

If there is a function $h\in\mathcal{H}$ such that
\begin{equation}\label{eq:multivar}
 g(x_1,x_2,\ldots,x_r) \; = \; 
  \int_{[0,1]} \prod_{i=1}^r h(x_i,y) \: d \mu(y) \: ,
\end{equation}
then $t(G,g) = t({\rm Inc}(G),h)$. 

Similarly to \cref{th:sub}, 
if ${\rm Inc}(G)$ satisfies \cref{conj:SC}, 
then \cref{eq:CIr} holds for functions $g$ 
that have representation \cref{eq:multivar} 
with nonnegative $h$. 
Similarly to \cref{th:norming}, 
if ${\rm Inc}(G)$ is norming (it requires $r$ to be even), 
then \cref{eq:CIr} holds for functions $g$ 
that have representation \cref{eq:multivar}, 
where $h\in\mathcal{H}$ can take negative values.

\subsection*{Acknowledgements}
The author would like to thank Joonkyung Lee 
for his valuable comments and suggestions. 

\subsection*{Acknowledgements}

\appendix

\section{Proofs of \cref{th:multipath,th:cycle}}

\begin{proof}[\bf{Proof of \cref{th:multipath}}]
The proof would be easy if all $k_i$ were even: 
then we could use the Cauchy--Schwarz inequality 
and the fact that the tree formed by paths of length 
$k_1/2,k_2/2,\ldots,k_r/2$ connected at one endpoint 
is an extra-good graph. 
To deal with odd values of $k_i$, 
we are going to subdivide the middle edge of the path. 
We assume that $k_1,\ldots,k_s$ are odd, 
and $k_{s+1},\ldots,k_r$ are even ($0 \leq s\leq r$). 
We denote the vertices of the $i$th path by 
$v_{i0},v_{i1},\ldots,v_{ik_i}$.
Gluing together vertices $v_{i0}$ for all $i$ 
produces vertex $v_0$. 
Gluing together vertices $v_{ik_i}$ 
produces vertex $v_\infty$. 
The total number of vertices is 
$n = 2 + \sum_{i=1}^r (k_i - 1)$, 
and the number of edges is $m = \sum_{i=1}^r k_i$.

Let $g$ be a doubly negative function: 
$g(x,y) = \int_{[0,1]} h(x,z)h(y,z) d\mu(z)$. 
Let $f_0,f_\infty,f_{ij}$ 
($1 \leq i \leq r$, $1 \leq j \leq k_i-1$) 
be bounded measurable nonnegative functions on $[0,1]$. 
Let
\[
  f \; = \; \left( 
    f_0 \cdot f_\infty \cdot
    \prod_{i=1}^r \prod_{j=1}^{k_i-1} f_{ij}
  \right)^{1/(2m)}
\]
For $i=1,2,\ldots,r$, denote
\[
      P_i \; = \; 
      g(x_0,x_{i1}) 
      \left(
      \prod_{j=2}^{k_i-1} g(x_{i,j-1},x_{i,j}) 
      \right)
      g(x_{i,k_i-1},x_\infty)
      \prod_{j=1}^{k_i-1} f_{i,j}(x_{i,j}) 
      \; .
\]
We need to prove $I \geq J^m$, where  
\[
  I =
  \int_{[0,1]^n} f_0(x_0) f_\infty(x_\infty) \prod_{i=1}^r P_i
    \: d\mu^n ,
    \;\;\;\;\;\;
  J =  
  \int_{[0,1]^2} f(x) g(x,y) f(y) \: d\mu^2 .
\]
We can assume that $f_0=f_\infty$, and 
$f_{i,j} = f_{i,k_i-j}$ for all $i$ and $j$. 
Indeed, we could swap $f_0$ with $f_\infty$, and 
$f_{i,j}$ with $f_{i,k_i-j}$ for all $i,j$ 
(this would not change the value of $I$), 
and then apply the Cauchy--Schwarz inequality 
to the geometric mean of both expressions for $I$.

For $i=1,2,\ldots,s$ (that is when $k_i$ is odd), 
we can add a variable $z_i$ and replace 
$g(x_{i,(k_i-1)/2},\:
   x_{i,(k_i+1)/2})
$ 
in the expression for $P_i$ with the product 
of $h(x_{i,(k_i-1)/2},z_i)$ and $h(x_{i,(k_i+1)/2},z_i)\:$:
\begin{multline*}
  P_i \; = \; 
      g(x_0,x_{i,1}) 
      \left(
      \prod_{j=2}^{(k_i-1)/2} g(x_{i,j-1},x_{i,j})
      \right)
      h(x_{i,(k_i-1)/2},z_i) \:
      h(x_{i,(k_i+1)/2},z_i) \\ 
      \times
      \left(
      \prod_{j=(k_i+3)/2}^{k_i-1} g(x_{i,j-1},x_{i,j}) 
      \right)
      g(x_{i,k_i-1},x_\infty)
      \prod_{j=1}^{k_i-1} f_{i,j}(x_{i,j})
      \; .
\end{multline*}
Now in the expression for $I$ we have to integrate 
over $n+s$ variables. 
Define for $i \leq s$,
\[
  Q_i \; = \;
      g(x_0,x_{i1}) 
      \prod_{j=2}^{(k_i-1)/2} g(x_{i,j-1},x_{i,j}) \;
      h(x_{i,(k_i-1)/2},z_i)
      \prod_{j=1}^{(k_i-1)/2} f_{i,j}(x_{i,j})
      \: ,
\]
and for $i > s$,
\[
  Q_i \; = \;
      g(x_0,x_{i,1}) \:
      \prod_{j=2}^{k_i/2} g(x_{i,j-1},x_{i,j})
      \prod_{j=1}^{k_i/2-1} f_{i,j}(x_{i,j}) \;
      \left(f_{i,k_i/2}(x_{i,k_i/2})\right)^{1/2}
      .
\]
Let $A(z_1,\ldots,z_s,\: x_{s+1,k_{s+1}/2},\ldots,x_{r,k_r/2})$ 
denote the integral of 
$f_0(x_0) \prod_{i=1}^r Q_i$ 
over all variables except 
$z_1,\ldots,z_s,\: x_{s+1,k_{s+1}/2},\ldots,x_{r,k_r/2}$. 
Set
\begin{multline*}
  B(y_1,\ldots,y_s,\: x_{s+1,k_{s+1}/2},\ldots,x_{r,k_r/2})
    \\ = \;
  \int_{[0,1]^s} \!
  A(z_1,\ldots,z_s,\: x_{s+1,k_{s+1}/2},\ldots,x_{r,k_r/2})
  \prod_{i=1}^s h(y_i,z_i) \: d\mu(z_1) \cdots d\mu(z_s)
    \: .
\end{multline*}
Then 
\[
  I \; = \;
  \int_{[0,1]^r}
    A(z_1,\ldots,z_s,\: x_{s+1,k_{s+1}/2},\ldots,x_{r,k_r/2})^2
    d\mu^r \: ,
\]
and
\begin{align*}
  I \cdot J^s \; & = \;
  I \cdot \int_{[0,1]^s} \prod_{i=1}^s
    \left( \int_{[0,1]} f(y_i) \: h(y_i,z_i) \: d\mu(y_i)
    \right)^2 d\mu(z_1) \cdots d\mu(z_s)
  \\ & \geq \;
  \left( \int_{[0,1]^r}
    B(y_1,\ldots,y_s,\: x_{s+1,k_{s+1}/2},\ldots,x_{r,k_r/2})
    \: d\mu^r
  \right)^2 .
\end{align*}
Notice that $\int_{[0,1]^r} B d\mu^r$ 
is the left hand side of \cref{eq:extra-good} 
for the tree $G$ formed by 
paths of length 
$(k_1+1)/2,\ldots,(k_s+1)/2,\: k_{s+1}/2,\ldots,k_r/2$ 
connected at one endpoint. 
By \cref{th:unicyclic}, $G$ is extra-good, so 
$\int_{[0,1]^r} B d\mu^r \geq J^{{\rm e}(G)}$. 
As ${\rm e}(G) = (m+s)/2$, we get
\[
  I \; \geq \; J^{-s} \left(\int_{[0,1]^r} B d\mu^r\right)^2
  \; \geq \; J^m \; .
\]
\end{proof}

\begin{proof}[\bf{Proof of \cref{th:cycle}}]
Let $G_0$ be formed by attaching a $k$-edge cycle 
$(v_0,v_1,\ldots,v_{k-1},\linebreak v_k=v_0)$ 
to a vertex $v$ of graph $G$, so $v=v_0=v_k$.
Let $g$ be a doubly nonnegative function. 
There exists $h \in {\mathcal H}$ such that
$g(x,y) = \int_{[0,1]} h(x,z)h(y,z) d\mu(z)$. 
Assign bounded measurable nonnegative functions on $[0,1]$
to all vertices of $G$, 
and denote by $F(x)$ the integral on the left hand side of 
\cref{eq:extra-good} taken over all variables except the one 
that corresponds to $v$. 
Assign bounded measurable nonnegative functions 
$f_1,f_2,\ldots,f_{k-1}$ to vertices $v_1,v_2,\ldots,v_{k-1}$. 
For functions $\gamma_1,\gamma_2,\ldots,\gamma_r$, denote 
\[
  I_r(x_0,x_r; \gamma_1,\gamma_2,\ldots,\gamma_r) \: =
  \int_{[0,1]^{r-1}} \prod_{i=1}^r \left(
    g(x_{i-1},x_i) \: \gamma_i(x_i)
  \right) d\mu(x_1) \cdots d\mu(x_{r-1}) .
\]
Let $I_{G_0}$ be the value of the integral 
on the left hand side of \cref{eq:extra-good} for $G_0$. 
Then
\[
  I_{G_0} \: =
  \int_{[0,1]} F(x) \: I_k(x,x; f_1,f_2,\ldots,f_{k-1},1)
  \: d\mu(x) \: .
\]
Let $f_0$ be the product of all ${\rm v}(G_0)$ functions 
assigned to vertices of $G_0$. 
We need to prove $I_{G_0} \geq J^{{\rm e}(G_0)}$, 
where $J=\int_{[0,1]^2} f(x)g(x,y)f(y)d\mu^2$ 
and $f=(f_0)^{1/(2{\rm e}(G_0))}$. 

Set $\gamma_i = \sqrt{f_i f_{k-i}}$, 
so $\gamma_{k-i} = \gamma_i$. 
By the Cauchy--Schwarz inequality, 
\begin{align*}
  I_k(x,x; f_1,f_2,\ldots,f_{k-1},1)^2
  \: & = \:
  I_k(x,x; f_1,f_2,\ldots,f_{k-1},1) \cdot
  I_k(x,x; f_{k-1},f_{k-2},\ldots,f_1,1)
  \\ & \geq \:
  I_k(x,x; \gamma_1,\gamma_2,\ldots,\gamma_{k-1},1)^2 .
\end{align*}
If $k=2a$, then
\begin{eqnarray*}
  I_{2a}(x,x; \gamma_1,\ldots,\gamma_{2a-1},1) & = &
  \int_{[0,1]} I_a(x,y;
    \gamma_1,\ldots,\gamma_{a-1},\sqrt{\gamma_a})^2
  \: d\mu(y) 
  \\ & \geq &
  \left(
  \int_{[0,1]} I_a(x,y;
    \gamma_1,\ldots,\gamma_{a-1},\sqrt{\gamma_a})
  \: d\mu(y) 
  \right)^2 .
\end{eqnarray*}
Construct graph $G_a$ from graph $G$ 
by attaching two disjoint $a$-edge paths 
to vertex $v$. 
By \cref{th:leaf}, $G_a$ is extra-good. 
Assign functions 
$\gamma_1,\ldots,\gamma_{a-1},
\sqrt{\gamma_a}$ 
to vertices of each of the two paths. 
Let $I_{G_a}$ be the value of the integral 
on the left hand side of \cref{eq:extra-good} for $G_a$. 
Then
\[
  I_{G_a} \: =
  \int_{[0,1]} F(x) \left(
  \int_{[0,1]} I_a(x,y;
    \gamma_1,\ldots,\gamma_{a-1},\sqrt{\gamma_a}) \: d\mu(y)
  \right)^2 d\mu(x)
  \: \leq \: I_{G_0} \: .
\]
The right hand side of \cref{eq:extra-good} 
is the same for both $G_0$ and $G_a$, 
so we are done with the case $k=2a$. 
If $k=2a+1$, we have
\begin{multline*}
  I_{2a+1}(x,x; \gamma_1,\ldots,\gamma_{2a},1) 
  \\ = 
  \int_{[0,1]^2} I_a(x,y'; \gamma_1,\ldots,\gamma_a) \:
                 g(y',y'') \:
                 I_a(x,y''; \gamma_{2a},\ldots,\gamma_{a+1})
  \: d\mu(y') d\mu(y'')
  \\ =
  \int_{[0,1]} \left(
    I_a(x,y; \gamma_1,\ldots,\gamma_a) \: h(y,z) \: d\mu(y)
  \right)^2 d\mu(z) \: ,
\end{multline*}
and
\[
  J \: = \int_{[0,1]^2} f(s') \: g(s',s'') \: f(s'') \: d\mu^2
  = \int_{[0,1]} \left( \int_{[0,1]} h(s,z) \: f(s) \: d\mu(s)
  \right)^2 d\mu(z) \: .
\]
Hence, by the Cauchy--Schwarz inequality,
\begin{multline*}
  I_{2a+1}(x,x; \gamma_1,\ldots,\gamma_{2a},1) \cdot J \: \geq
  \\ \geq \left( \int_{[0,1]^3}
    I_a(x,y; \gamma_1,\ldots,\gamma_a) 
    \: h(y,z) \: h(s,z) \: f(s) \: d\mu(y) d\mu(z) d\mu(s)
  \right)^2
  \\ = \left(
    \int_{[0,1]} I_{a+1}(x,s; \gamma_1,\ldots,\gamma_a,f) 
      \: d\mu(s) 
    \right)^2 .
\end{multline*}
Construct $G_{a+1}$ from $G$ 
by attaching two disjoint $(a+1)$-edge paths 
to vertex $v$. 
By \cref{th:leaf}, $G_{a+1}$ is extra-good. 
Assign functions $\gamma_1,\ldots,\gamma_a,f$ 
to vertices of each of the two paths. 
Let $I_{G_{a+1}}$ be the value of the integral 
on the left hand side of \cref{eq:extra-good} for $G_{a+1}$. 
Then
\[
  I_{G_{a+1}} \: =
  \int_{[0,1]} F(x) \left(
  \int_{[0,1]} I_{a+1}(x,s;
    \gamma_1,\ldots,\gamma_a,f) \: d\mu(s)
  \right)^2 d\mu(x)
  \: \leq \: I_{G_0} \cdot J \: .
\]
The right hand side of \cref{eq:extra-good} 
is equal to $J^{{\rm e}(G_0)}$ for $G_0$, 
and $J^{{\rm e}(G_{a+1})} = J \cdot J^{{\rm e}(G_0)}$ 
for $G_{a+1}$. 
As \cref{eq:extra-good} holds for $G_{a+1}$,
it holds for $G_0$, too. 
\end{proof}


\begin{thebibliography}{99}

\bibitem{Adams:2008}
P. Adams, D. Bryant, and M. Buchanan. 
\newblock A survey on the existence of {$G$}-designs. 
\newblock \emph{J. Combin. Designs}, 16(5):373--410, 2008. 
\doi{10.1002/jcd.20170}.

\bibitem{Atkinson:1960}
F. V. Atkinson, G. A. Watterson, and P. A. D. Moran.
\newblock A matrix inequality.
\newblock \emph{Quarterly J. of Math.}, 
11(42):137--140, 1960.
\doi{10.1093/qmath/11.1.137}.

\bibitem{Berman:2003} 
A. Berman and N. Shaked-Monderer. \newblock 
\newblock \emph{Completely Positive Matrices}, World Scientific, 2003. 
\doi{10.1142/5273}.

\bibitem{Blakley:1965}
G. R. Blakley and P. Roy. 
\newblock H\"{o}lder type inequality for symmetric matrices with
nonnegative entries. 
\newblock \emph{Proc. Amer. Math. Soc.}, 16(6):1244--1245, 1965.
\doi{10.1090/S0002-9939-1965-0184950-9}.

\bibitem{Blinco:2004}
A. Blinco. 
\newblock Theta graphs, graph decompositions and related graph labelling techniques. 
\newblock \emph{Bull. Austral. Math. Soc.}, 69(1):173--175, 2004. 
\doi{10.1017/S0004972700034377}.

\bibitem{Conlon:2010} 
D. Conlon, J. Fox, and B. Sudakov. 
\newblock An approximate version of Sidorenko's conjecture. 
\newblock \emph{Geom. Funct. Anal.}, 20(6):1354--1366, 2010. 
\doi{10.1007/s00039-010-0097-0}.

\bibitem{Conlon:2018}
D. Conlon, J. H. Kim, C. Lee, and J. Lee. 
\newblock Some advances on Sidorenko's conjecture. 
\newblock \emph{J. London Math. Soc.}, 98(3):593--608, 2018.  \doi{10.1112/jlms.12142}.

\bibitem{Conlon:2017}
D. Conlon and J. Lee. 
\newblock Finite reflection groups and graph norms. 
\newblock \emph{Adv. Math.}, 315:130--165, 2017.  \doi{10.1016/j.aim.2017.05.009}.

\bibitem{Conlon:2019}
D. Conlon and J. Lee. 
\newblock Sidorenko's conjecture for blow-ups. 
\newblock \emph{Discrete Analysis, to appear}.  
\newblock \arxiv{1809.01259}, 2018.

\bibitem{Hatami:2010}
H. Hatami.
\newblock Graph norms and Sidorenko's conjecture. 
\newblock \emph{Israel J. Math.}, 175(1):125--150, 2010. 
\doi{10.1007/s11856-010-0005-1}.

\bibitem{Jagger:1996}
C. Jagger, P. \v{S}\v{t}ovi\v{c}ek, and A. Thomason. 
\newblock Multiplicities of subgraphs. 
\newblock \emph{Combinatorica}, 16(1):123--141, 1996. 
\doi{10.1007/BF01300130}.

\bibitem{Kim:2016}
J. H. Kim, C. Lee, and J. Lee. 
\newblock Two approaches to Sidorenko's conjecture. 
\newblock \emph{Trans. Amer. Math. Soc.}, 368(7):5057--5074, 2016. 
\doi{10.1090/tran/6487}.

\bibitem{Kohayakawa:2010}
Y. Kohayakawa, B. Nagle, V. R\"{o}dl, and M. Schacht. 
\newblock Weak hypergraph regularity and linear hypergraphs. 
\newblock \emph{J. Combin. Theory Ser. B}, 100(2):151--160, 2010.
\doi{10.1016/j.jctb.2009.05.005}.

\bibitem{Lee:2019b}
J. Lee.
\newblock On some graph densities in locally dense graphs.
\newblock \emph{Random Struct Alg.}, 58:322--344, 2021.
\doi{10.1002/rsa.20974}.

\bibitem{Li:2011}
J. L. X. Li and B. Szegedy. 
\newblock On the logarithmic calculus and Sidorenko's conjecture. 
\newblock \arxiv{1107.1153}, 2011.

\bibitem{Lovasz:2010}
L. Lov\'{a}sz. 
\newblock \emph{Large networks and graph limits}, 
volume 60 of \emph{Colloquium Publications}. AMS, 2012. 
\doi{10.1090/coll/060}.

\bibitem{Lovasz:2011}
L. Lov\'{a}sz. 
\newblock Subgraph densities in signed graphons and the local Sidorenko conjecture. 
\newblock \emph{Electr. J. Combin.}, 18(1) \#P127, 2011.

\bibitem{Mulholland:1959}
H. P. Mulholland and C. A. B. Smith. 
\newblock An inequality arising in genetical theory. 
\newblock \emph{Amer. Math. Monthly}, 66(8):673--683, 1959. 
\doi{10.2307/2309342}.

\bibitem{Parczyk:2014}
O. Parczyk. 
\newblock On Sidorenko's conjecture, Master’s thesis, 
Freie Universit\"{a}t Berlin, 2014. 
\url{http://www.uni-frankfurt.de/58522166}.

\bibitem{Reiher:2014}
C. Reiher. 
\newblock Counting odd cycles in locally dense graphs. 
\newblock \emph{J. Combin. Theory Ser. B}, 105:1--5, 2014.
\doi{10.1016/j.jctb.2013.12.002}.

\bibitem{Sidorenko:1992}
A. Sidorenko. 
\newblock Inequalities for functionals generated by bipartite graphs. 
\newblock \emph{Discrete Math. Appl.}, 2(5):489--504, 1992. 
\doi{10.1515/dma.1992.2.5.489}.

\bibitem{Sidorenko:1993}
A. Sidorenko. 
\newblock A correlation inequality for bipartite graphs. 
\newblock \emph{Graphs and Combinatorics}, 9(2):201--204, 1993. 
\doi{10.1007/BF02988307}.

\bibitem{Sidorenko:1991}
A. Sidorenko. 
\newblock An analytic approach to extremal problems for graphs and hypergraphs. \newblock In 
\emph{Extremal Problems for Finite Sets (Visegr\'{a}d, 1991)}, 
volume 3 of \emph{Bolyai Soc. Math. Stud.}, pages 423--455. 
J\'{a}nos Bolyai Math. Soc., Budapest, 1994.

\bibitem{Szegedy:2014}
B. Szegedy. 
\newblock An information theoretic approach to Sidorenko's conjecture.
\newblock \arxiv{1406.6738}, 2014.

\bibitem{Szegedy:2015}
B. Szegedy. 
\newblock Sparse graph limits, entropy maximization and transitive graphs.
\newblock \arxiv{1504.00858}, 2015.

\end{thebibliography}
\end{document}